\documentclass[11pt,reqno]{amsart}
\usepackage[foot]{amsaddr}

\usepackage{hyperref,url}

 \usepackage{float}
 \usepackage{mathtools}
\usepackage{amssymb}
\usepackage{bbm}
\usepackage{amssymb}
\usepackage{amsmath, amsfonts}
\usepackage{amscd}
\usepackage{amsthm}
\usepackage{setspace}
\usepackage{enumerate}
\usepackage{graphicx}
\usepackage{mathtools}
\usepackage{tcolorbox} 
\usepackage{subfigure} 
\usepackage{siunitx}
\usepackage{tikz-cd}
\usepackage{bm}
\numberwithin{equation}{section}

\usepackage{mathtools}
\usepackage[tableposition=top]{caption}
\usepackage{booktabs,dcolumn}
\usepackage{bigints}

\renewcommand\H{\mathbb{H}}

\newcommand\T{\mathbb{T}}

\renewcommand\k{\mathbf{k}}

\newcommand\A{\mathbf{A}}

\renewcommand\hat{\widehat}

\newcommand\eps{\varepsilon}

\DeclareMathOperator{\supp}{supp}

\DeclareMathOperator{\sech}{sech}

\DeclareMathOperator{\rad}{rad}
\DeclareMathOperator{\Span}{Span}

\newcommand{\quash}[1]{}

\renewcommand{\H}{\textbf H}

\def\XXint#1#2#3{{\setbox0=\hbox{$#1{#2#3}{\int}$}
     \vcenter{\hbox{$#2#3$}}\kern-.5\wd0}}

\title[Fourier inequalities and sign uncertainty II]{Fourier inequalities and sign uncertainty II}
\author{Roni Edwin}
\newtheorem{theorem}{Theorem}[section]

\newtheorem{lemma}[theorem]{Lemma}
\newtheorem{problem}[theorem]{Problem}

\begin{document}

\address{Yale University, New Haven, Connecticut 06511}
\email{roni.edwin@yale.edu}


\begin{abstract}
Building on our earlier work, we present an $L^2$-$L^1$ inequality for radial Fourier eigenfunctions, relating this to the sign uncertainty principle of  Bourgain, Clozel and Kahane. We also propose a generalisation of the aforementioned uncertainty principle, and draw a connection with entropic uncertainty.
\end{abstract}

\maketitle

\tableofcontents

\renewcommand{\theenumi}{\alph{enumi}}
\section{Introduction}
Recently, OpenAI released a preprint \cite{OpenAI2026TenAdvances} in which they resolved a number of open problems in mathematics. Most relevant to us is their resolution of the asymptotic behaviour of the optimal constant in the sign uncertainty principle (\cite{bourgain2010principe}, \cite{cohn2019optimal}, \cite{gonccalves2017hermite}). If for $\sigma \in \{-1,1\}$, $\A_\sigma(d)$ denotes the best constant in the $\sigma$-sign uncertainty principle, they showed that\begin{align}
    \lim_{d \to \infty}\frac{\mathbf{A}_+(d)}{\sqrt{d}}=\lim_{d \to \infty}\frac{\mathbf{A}_-(d)}{\sqrt{d}}=\frac1{\pi},
    \label{person}
\end{align} matching the conjectured asymptotic in \cite{modbootstrap}. Regarding the progress towards this result, we note that early work \cite{bourgain2010principe}  showed that \begin{align}
    \frac1{\sqrt{2\pi e}}\le \liminf_{d \to \infty} \frac{\mathbf{A}_+(d)}{\sqrt{d}}\le \limsup_{d \to \infty}\frac{\mathbf{A}_+(d)}{\sqrt{d}}\le \frac1{\sqrt{2\pi }},
    \label{worka}
\end{align} and for a while these seemed to be the best asymptotic bounds, though it was conjectured in \cite{modbootstrap} that $\lim_{d \to \infty}\frac{\A_\pm(d)}{\sqrt{d}}=\frac1{\pi}$. We did manage to improve the lower bound to \begin{align}
    \liminf_{d \to \infty}\frac{\A_\pm(d)}{\sqrt{d}}\ge \frac1{2\sqrt{\pi}}
    \label{feeling}
\end{align} in \cite{EdwFourier}, with \eqref{person} finally being proven in \cite{OpenAI2026TenAdvances}.
 This paper sharpens our earlier result, exploring the methods used to obtain \eqref{person}, and proposes a generalisation of the sign uncertainty principle.  Before diving in, some background for the uninitiated:

The sign uncertainty principle, originally posed in  \cite{bourgain2010principe}, with a complementary $-1$-uncertainty principle developed in \cite{cohn2019optimal}, is as follows. On $\mathbb{R}^d$, we let $\mathcal{A}_+(d)$ be the set of non-zero functions $f\colon \mathbb{R}^d \to \mathbb{R}$ such that $f,\hat{f} \in L^1(\mathbb{R}^d;\mathbb{R})$ (using continuous representatives), with $f$ being eventually non-negative and $f(0)\le 0$, and similarly, $\hat{f}$ is eventually non-negative with $\hat{f}(0)\le 0$.

One notes membership in $\mathcal{A}_+(d)$ implies neither $f$ nor $\widehat{f}
$ can be non-negative outside an arbitrarily small neighbourhood of the origin. That is, the following quantities \begin{align*}
    A(f)&\coloneqq \inf\left\{r>0: f(x)\ge 0 \text{ for all }|x|\ge r\right\}, \\
    A(\widehat{f})&\coloneqq \inf\left\{r>0: \widehat{f}(x)\ge 0 \text{ for all }|x|\ge r\right\}
\end{align*} are strictly positive. A natural quantity to study is the product \begin{align}
\sqrt{A(f)A(\widehat{f})},
    \label{prod}
\end{align} as it is invariant under dilation $x \to \lambda x$. In particular, one can choose an appropriate dilation and consider $f+ \hat{f}$, in which case one restricts the consideration of the product $A(f)A(\hat{f})$ to $ 1$-eigenfunctions. Moreover, one can impose the condition $f(0)=0$, since the function $x\mapsto f(x)+(-f(0))e^{-\pi |x|^2}$ is still a $+1$-eigenfunction, with $A\mathopen{}(f+(-f(0))e^{-\pi |\cdot|^2})\mathclose{}\le A(f)$. This leads to the following $+1$-eigenfunction uncertainty principle on $\mathbb{R}^d$:
\begin{problem}[$+1$-uncertainty principle]
    Minimise $A(g)$ over all $g\colon \mathbb{R}^d \to \mathbb{R}$ such that \begin{enumerate}
        \item $g \in L^1\mathopen{}\left(\mathbb{R}^d\right)\mathclose{}\mathopen{}\setminus\mathclose{} \{\mathbf{0}\}$ and $\widehat{g}= g$,
        \item $g(0)=0$ and $g$ is eventually non-negative.
    \end{enumerate}
    \label{dancefloor}
\end{problem} Call the minimal value of $A(g)$ in this problem $\mathbf{A}_+(d)$.
Theorem $3$ in \cite{gonccalves2017hermite} shows $\A_+(d)$ can be realised by a non-zero integrable $+1$-eigenfunction. Thus Problem \ref{dancefloor} is equivalent to minimising the product \eqref{prod}. We note this has been solved explicitly in $d=12$ in \cite{cohn2019optimal}, with $\A_+(12)=\sqrt{2}$. Lee \cite{lee2026positivequasimodularformssign} adapted their methods, obtaining the upper bound \begin{align*}
    \A_+(d)\le \sqrt{2\left\lfloor \frac{d}{16}\right\rfloor+2},
\end{align*}
for dimensions $d$ divisible by $4$. 

The $-1$-uncertainty principle is defined analogously:
\begin{problem}[$-1$-uncertainty principle]
    Minimise $A(g)$ over all $g\colon \mathbb{R}^d \to \mathbb{R}$ such that \begin{enumerate}
        \item \label{i} $g \in L^1\mathopen{}\left(\mathbb{R}^d\right)\mathclose{}\mathopen{}\setminus\mathclose{} \{\mathbf{0}\}$ and $\widehat{g}=-g$,
        \item \label{ii} $g(0)=0$ and $g$ is eventually non-negative.
    \end{enumerate}
    \label{dancefloor2}
\end{problem} 
We similarly call the minimal value of $A(g)$ in the above problem $\mathbf{A}_-(d)$, so \begin{align*}
    \mathbf{A}_-(d)\coloneqq \inf\left\{A(g): g\text{ satisfies \ref{i} and \ref{ii}}\right\}.
\end{align*}
This problem has been solved for $d=1,8,24$, as a consequence of the linear programming bounds (Theorem 3.2 in \cite{cohn2003new}) for sphere packing (see \cite{cohn2019optimal} for further discussion).

We now explore how the lower bounds for $\A_\sigma(d)$, $\sigma =+,-$, were developed. As a simplification, we note we can restrict our attention to radial functions. If $\mathcal{R}_d$ denotes the radialisation operator on $\mathbb{R}^d$, so \begin{align}
    \mathcal{R}_dg(x)=\int_{O(d)}g(Ux)\nu_d(\textup{d}U),
    \label{radialisationdef}
\end{align} where $\nu_d$ is the normalised Haar measure on $O(d)$, replacing $f$ with $\mathcal{R}_df$ does not increase either of $A(f)$, $A(\pm \hat{f})$. Moreover, $\mathcal{R}_df\neq 0$, since otherwise $f$ and $\hat{f}$ would both be compactly supported.
One starts by noting that given a $\sigma$-eigenfunction $g$ as in Problems \ref{dancefloor} and \ref{dancefloor2}, if $A(g)$ denotes the radius of its last sign change, then (assuming $g$ eventually becomes non-negative), this forces the negative part of $g$ to lie inside the ball of radius $A(g)$; That is, writing $g=g_+-g_-$, $g_\sigma=\max(\sigma g, 0)$, one necessarily has $\supp(g_-) \subset \overline{B_{A(g)}}$, thus \begin{align}
    \frac1{2}\lVert g\rVert_{L^1}=\int_{B_{A(g)}}g_-\le \int_{B_{A(g)}}|g|.
    \label{lasse}
\end{align} This inequality suggests that a useful approach for obtaining lower bounds for $\A_\sigma(d)$ is controlling this `mass concentration', \begin{align}
    \sup_{\substack{\hat{g}=\sigma g \text{ radial } \\ g(0)=0,\ g \neq 0}}\frac1{\lVert g\rVert_{L^1}}\int_{B_r}|g|.
    \label{baby}
\end{align}
 Indeed, if for some $\tau>0$, one has \begin{align*}
    \frac1{\lVert g\rVert_{L^1}}\int_{B_{\tau\sqrt{d}}}|g|=o(1),
\end{align*} uniformly for all radial $\pm 1$-eigenfunctions $g$ vanishing at the origin as $d \to \infty$,  combining this bound with the inequality \eqref{lasse}  implies \begin{align*}
    \liminf_{d \to \infty}\frac{\mathbf{A}_\sigma(d)}{\sqrt{d}}\ge \tau.
\end{align*} This approach has proven useful in lower bounding the limit inferior above. The initial lower bound in \eqref{worka} was obtained via this approach. By applying H\"{o}lder's inequality one obtains for an applicable eigenfunction $g$, \begin{align*}
    \int_{B_r}|g|\le |B_r|\lVert g\rVert_\infty\le |B_r|\lVert g\rVert_{L^1},
\end{align*} thus if $|B_r|\le 1/2$, this forces most of the mass of $g$ outside the ball. Asymptotically, this leads to the lower bound in \eqref{worka}, \begin{align*}
    \liminf_{d \to \infty}\frac{\A_\sigma(d)}{\sqrt{d}}\ge \frac1{\sqrt{2\pi e}}.
\end{align*} Our initial improvement on this lower bound, \eqref{feeling},  came from applying Cauchy-Schwarz to \eqref{baby}, obtaining \begin{align}
   \mathopen{}\left( \int_{B_r}|g|\right)^2\mathclose{}\le |B_r|\lVert g\rVert_{L^2}^2,
   \label{asabut}
\end{align} and bounding the $L^2$-norm by the $L^1$-norm. This was the focus of our Theorem 1.5: \begin{theorem}[Theorem 1.5 in \cite{EdwFourier}]
    Let $\mathcal{H}^d$ be the subspace of $L^1\mathopen{}\left(\mathbb{R}^d\right)$ defined by \begin{align}
\mathcal{H}^d\coloneqq \left\{f \in L^1\mathopen{}\left(\mathbb{R}^d;\mathbb{R}\right)\mathclose{}: \widehat{f} \in L^1\mathopen{}\left(\mathbb{R}^d;\mathbb{R}\right)\mathclose{}\right\}.
    \label{mathcalHdef}
\end{align}
Then
\begin{align}
        \left\lVert f\right\rVert_{L^2}^2\le \left(\frac2{e}\right)^{\frac{d}{2}}\left\lVert f\right\rVert_{L^1} \big\lVert \widehat{f}\big\rVert_{L^1},
        \label{mistakes}
    \end{align} for all $f \in \mathcal{H}^d$.
\label{howtolove}
\end{theorem} This inequality combined with \eqref{asabut} leads to the lower bound in \eqref{feeling}. Unfortunately, this method is not sharp, in the sense that the best constant in \eqref{mistakes} does not produce the bound in \eqref{person}. That is, if \begin{align}
    \gamma\coloneqq \limsup_{d\to \infty}\mathopen{}\left(\sup_{f\in \mathcal{H}_{\rad}^d\mathopen{}\setminus\{\mathbf{0}\}}\frac{\lVert f\rVert_{L^2}^2}{\lVert f\rVert_{L^1}\lVert \hat{f}\rVert_{L^1}}\right)^{\frac1{d}},
    \label{remaa}
\end{align} applying this to \eqref{asabut} leads to \begin{align*}
    \left(\int_{B_{\tau\sqrt{d}}}|g|\right)^2\le \left|B_{\tau\sqrt{d}}\right|(\gamma+o(1))^d\lVert g\rVert_{L^1}^2 \sim \frac{( \tau(\gamma+o(1))\sqrt{2\pi e})^d}{\sqrt{\pi d}}\lVert g\rVert_{L^1}^2,
\end{align*} from Stirling's approximation. This implies a lower bound of \begin{align}
    \liminf_{d \to \infty}\frac{\mathbf{A}_\sigma(d)}{\sqrt{d}}\ge \frac1{\gamma\sqrt{2\pi e}}.
    \label{tension}
\end{align} Numerically one has $\gamma \approx  0.789.. $, thus leading to a constant on the right-hand side of \eqref{tension} looking like \begin{align*}
    0.30654...<\frac1{\pi}\approx 0.318309886184.
\end{align*} So while it is not terribly inefficient, it does not recover \eqref{person}. Motivated by the result in \eqref{person}, we considered two questions, the first being, what is the optimal constant $\gamma$ in \eqref{remaa}? By adapting some of the ideas in Chapter 1 of \cite{OpenAI2026TenAdvances}, we were able to show: 
\begin{theorem}
Let $\mathcal{H}_{\textup{rad}}^d$ denote the space of  radial functions in $\mathcal{H}^d$, and set \begin{align*}
    w_d^{\rad}\coloneqq \sup_{\substack{f \in \mathcal{H}_{\textup{rad}}^d\mathopen{}\setminus\mathclose{}\{\mathbf{0}\}}}\frac{\lVert f\rVert_{L^2}^2}{\lVert f\rVert_{L^1}\lVert \hat{f}\rVert_{L^1}}.
\end{align*}
Then there exist constants $C,c>0$ such that \begin{align*}
    cd^{-2}2^{\frac{d}{2}}e^{-\frac{2dG}{\pi}}\le w_d^{\rad}\le Cd^{\frac1{2}}2^{\frac{d}{2}}e^{-\frac{2dG}{\pi}},
\end{align*} where $G=\sum_{k=0}^\infty \frac{(-1)^k}{(2k+1)^2}$ is Catalan's constant. Thus
\begin{align}
     \lim_{d \to \infty}\mathopen{}\mathopen{}\left(w_d^{\rad}\right)^{\frac1{d}}\mathclose{}=\sqrt{2}e^{-\frac{2G}{\pi}}\approx 0.789347821.
     \label{trainingseason}
\end{align}
    \label{mainthm2}
\end{theorem}

Another question is as follows. We note the mass concentration approach outlined earlier bounds the measure of the set where $g$ is non-positive. Specifically, suppose we have a $\sigma$-eigenfunction $g$ with $g(0)=0$ that is eventually non-negative, and let $S=\{g\le 0\}$. Then one has 
\begin{align}
    \frac1{4}\lVert g\rVert_{L^1}^2=\left(\int_Sg\right)^2\le |S|\lVert g\rVert_{L^2}^2\le (2/e)^{\frac{d}{2}}|S|\lVert g\rVert_{L^1}^2,
    \label{shroom}
\end{align} thus $|S|\ge {(2/e)^{-\frac{d}{2}}}/{4}$. This led us to the following question:
Given $f \in \mathcal{H}^d$ such that $f$ is non-positive on a set of finite measure, let \begin{align*}
    M_d(f)\coloneqq \left\{x \in \mathbb{R}^d: f(x)\le 0\right\}.
\end{align*}
\begin{problem}[Definite-set uncertainty principle]
Let $f \in \mathcal{H}^d\setminus \{\mathbf{0}\}$ (so $f$ is a real-valued, even integrable function on $\mathbb{R}^d$ whose Fourier transform is also integrable) be a Fourier $\sigma$-eigenfunction with $f(0)=0$. How small can $|M_d(f)|$ be?
    \label{generalmeas}
\end{problem} Call the infimum $\mathcal{M}_\sigma(d)$, so \begin{align*}
    \mathcal{M}_\sigma(d)=\inf\mathopen{}\left\{|M_d(f)|: f\in \mathcal{H}^d\mathopen{}\setminus \mathclose{}\{\mathbf{0}\},\ \hat{f}=\sigma f, \ f(0)=0\right\}.
\end{align*}
Unlike the sign uncertainty problem, this problem does not seem to admit an immediate reduction to radial functions; the effect of radialising on $|M_d(f)|$ is unclear. Nonetheless, a useful avenue for obtaining bounds on $\mathcal{M}_\sigma(d)$ may be obtaining inequalities of the form \begin{align*}
    \mathopen{}\left(  \int_{X}|f|\right)\mathclose{}\mathopen{}\left(  \int_{\Sigma}\left|\hat{f}\right|\right)\mathclose{}\le C(X,\Sigma)\lVert f\rVert_{L^1}\lVert \hat{f}\rVert_{L^1}.
\end{align*}
The analogous problem involving the $L^2$-norm has been studied extensively; see for example \cite{jaming2007nazarov}. The corresponding $L^1$ theory appears to be less developed; generally speaking the literature seems to be lacking comparatively on quantitative uncertainty principles involving the $L^p$ norms, for $p\neq 2$. 

There are a number of reasons $\mathcal{M}_\sigma(d)$ and $M_d(f)$ may be of interest to us. It seems a natural stepping point from the classical sign uncertainty principle. Moreover, while the elimination of the radialisation reduction may render the resulting problem more difficult, progress in this setting may reveal previously unexplored connections between different areas of Fourier analysis.


So far, our best lower bound for $\mathcal{M}_\sigma(d)$ is that implied by Theorem \ref{howtolove}: The Cauchy-Schwarz argument in \eqref{shroom} gives \begin{align*}
    \mathcal{M}_\sigma (d)\ge \frac{w_d^{-1}}{4} \ge  \frac1{4}\mathopen{}\left(\frac{e}{2}\right)^{\frac{d}{2}},
\end{align*} where \begin{align*}
    w_d\coloneqq \sup_{f \in \mathcal{H}^d\mathopen{}\setminus \{\mathbf{0}\}}\frac{\lVert f\rVert_{L^2}^2}{\lVert f\rVert_{L^1}\lVert \hat{f}\rVert_{L^1}}.
\end{align*}
A simple tensorisation argument implies $w_{n+m}\ge w_nw_m$, thus the limit $\lim_{d \to \infty}(w_d)^{\frac1{d}}$ exists. We do not know if in a specific dimension $d$, the supremum can be attained by radial functions, though calculations in a few select dimensions $d$ seem to suggest $w_d=w_d^{\rad}$; we were unable to find better non-radial candidates. Perhaps a more accessible question is if \begin{align*}
    \lim_{d \to \infty}\mathopen{}\left( w_d\right)^{\frac1{d}}=\limsup_{d\to \infty}(w_d^{\rad})^{\frac1{d}}.
\end{align*} 

There is also a connection between this problem and entropy. One notes the Cauchy-Schwarz approach can be modified to $L^p$, $p \in (1,\infty)$: One has with $p=1+\alpha>1$, $f$ a $\sigma$-eigenfunction with $f(0)=0$, \begin{align*}
    \frac{\lVert f\rVert_{L^1}}{2}=\int_{M_d(f)}f_-\le \lVert f\rVert_{L^p}\left\lVert \mathbbm{1}_{M_d(f)}\right\rVert_{L^{p^*}}=\lVert f\rVert_{L^p}|M_d(f)|^{\frac1{p^*}},
\end{align*} so \begin{align*}
    |M_d(f)|\ge 2^{-{p^*}}\mathopen{}\left(\frac{\lVert f\rVert_{L^1}}{\lVert f\rVert_{L^p}}\right)^{\frac{p}{p-1}}.
\end{align*} Now normalise $f$ so that $\lVert f\rVert_{L^p}=1$. Then \begin{align*}
   \lVert f\rVert_{L^1}= \int_{\mathbb{R}^d}|f|^{1+\alpha}|f|^{-\alpha}=\int_{\mathbb{R}^d}|f|^{1+\alpha}\exp\mathopen{}\left(-\alpha\log|f|\right)\mathclose{}\ge \exp\mathopen{}\left(-(p-1) \int_{\mathbb{R}^d}|f|^{p}\log|f|\right),
\end{align*} by Jensen's inequality. Thus \begin{align*}
    |M_d(f)|\ge 2^{-\frac{p}{p-1}}\exp\mathopen{}\left(-p\int_{\mathbb{R}^d}|f|^p\log|f|\right)\mathclose{}.
\end{align*} Indeed, the $p=2$ case recovers the asymptotic lower bound \eqref{feeling}, following from the classic $L^2$ entropic uncertainty principle (\cite{folland1997uncertainty}, Theorem 5.3). We expect (in the case of the sign uncertainty principle at least) this method to be exhaustive as $p \to 1$. This suggests general entropic $L^p$ uncertainty principles of the form \begin{align*}
    \H\mathopen{}\left(|f|^p/\lVert f\rVert_{L^p}^p\right)\mathclose{}+\H\mathopen{}\left(\big\lvert \hat{f}\big\rvert^p/\lVert \hat{f}\rVert_{L^p}^p\right)\mathclose{}\ge C_{p,d},
\end{align*} might be useful. While the literature is fairly exhaustive in the $p=2$ case, it seems to be comparatively lacking for $p \in [1,2)$.

The paper is organised as follows. The proof of Theorem \ref{mainthm2} is split across Sections \ref{three} and \ref{lowerbound}. Section \ref{dualform} includes numerical upper bounds for $(w_d^{\rad})^{\frac1{d}}$, and the corresponding lower bound on $\A_\sigma(d)$.
We now present the proof of Theorem \ref{mainthm2}.

\subsection{Disclosure on AI use}
ChatGPT-6.1 Sol was also used in generating the code for computing the numerics in Section \ref{dualform}, as well as proofreading for typographical/mathematical errors.

\section{Preliminaries}\label{areaimprov}
We remark that $\mathcal{S}_{\rad}(\mathbb{R}^d)$ is dense in $\mathcal{H}_{\rad}^d$ under the norm \begin{align}
\lVert g\rVert_{\mathcal{H}^d_{\rad}}\coloneqq \lVert g\rVert_{L^1}+\lVert \hat{g}\rVert_{L^1},
    \label{Hradnorm}
\end{align} so we can restrict our attention to Schwartz functions.
For our analysis we would like to restrict to the case where the dimension $d$ is even. We start by establishing: \begin{lemma}
Let $w_d^{\rad}$ be as defined in Theorem \ref{mainthm2}. Then $w_d^{\rad}\le w_{d-1}^{\rad}$ for all $d\ge 2$.
    \label{firstlemma}
\end{lemma}
\begin{proof}
    Consider any non-zero function $f \in \mathcal{H}_{\rad}^d\cap \mathcal{S}(\mathbb{R}^d)$. Without loss of generality, let $\widehat{f}$ not be compactly supported ($f$ and $\widehat{f}$ cannot both be compactly supported, established in \cite{jaming2007nazarov}). For each $s \in \mathbb{R}$, define $H_s\colon \mathbb{R}^{d-1} \to \mathbb{R}$ by $H_s(\omega)=\widehat{f}(\omega,s)$. Note $H_s$ is radial, since $f$ is radial, and $H_s$ is non-zero, since $\widehat{f}$ is not compactly supported.  $H_s$ is integrable, from the decay of $\hat{f}$. To see that $\widehat{H_s}$ is also integrable, observe that \begin{align*}
     \int_{\mathbb{R}^{d-1}}\left(\int_{\mathbb{R}} f\mathopen{}\left(q_1,q_2\right)\mathclose{}e^{-2\pi i  sq_2 }\textup{d}q_2\right)e^{-2\pi i \left\langle q_1,\omega\right\rangle}\textup{d}q_1&=\int_{\mathbb{R}^{d-1}\times \mathbb{R}}f\left(q_1,q_2\right)e^{-2\pi i \left(\left\langle q_1,\omega\right\rangle+q_2s\right)}\textup{d}q_2\textup{d}q_1 \\
       &=\int_{\mathbb{R}^d}f(q)e^{-2\pi i \left\langle (q_1,q_2), (\omega,s)\right\rangle}\textup{d}q=\widehat{f}\left(\omega,s\right).
    \end{align*} This tells us $H_s$ is the Fourier transform of the integrable function $ \bigintsss_{\mathbb{R}}f(\cdot,y)e^{-2\pi i  sy}\textup{d}y$,
    so  \begin{align*}
        \widehat{H_s}(x)=\int_{\mathbb{R}}f(x,y)e^{-2\pi i sy}\textup{d}y.
    \end{align*} Having established the integrability of $H_s$ and $\hat{H_s}$, we have \begin{align*}
        w^{\rad}_{d-1}&\ge \frac{\int_{\mathbb{R}^{d-1}}\left|H_s\right|^2}{\left(\int_{\mathbb{R}^{d-1}}\left|H_s\right|\right)\left(\int_{\mathbb{R}^{d-1}}\left|\widehat{H_s}\right|\right)},
    \end{align*} and by clearing the denominator and using the definition of $H_s$ and the formula for $\widehat{H_s}$, \begin{align*}
        \int_{\mathbb{R}^{d-1}}\left|\widehat{f}(\omega,s)\right|^2\textup{d}\omega\le w_{d-1}^{\rad}\left(\int_{\mathbb{R}^{d-1}}\left|\widehat{f}(\omega,s)\right|\textup{d}\omega\right)\left(\int_{\mathbb{R}^{d-1}}\left|\int_{\mathbb{R}}f(x,y)e^{-2\pi i sy}\textup{d}y\right|\textup{d}x\right).
    \end{align*} Integrating over $s \in \mathbb{R}$ then gives
     \begin{align*}
        \int_{\mathbb{R}^d}\left|f\right|^2&\le w^{\rad}_{d-1} \bigintsss_{\mathbb{R}}\left(\int_{\mathbb{R}^{d-1}}\left|\widehat{f}(\omega,s)\right|\textup{d}\omega\right)\left(\int_{\mathbb{R}^{d-1}}\left|\int_{\mathbb{R}}f(x,y)e^{-2\pi i sy}\mathclose{}\textup{d}y\right|\textup{d}x\right)\mathclose{}\textup{d}s \\
        &\le w_{d-1}^{\rad}\bigintsss_{\mathbb{R}}\left(\int_{\mathbb{R}^{d-1}}\left|\widehat{f}(\omega,s)\right|\textup{d}\omega\right)\left(\int_{\mathbb{R}^{d-1}}\int_{\mathbb{R}}\left|f(x,y)\right|\mathclose{}\textup{d}y\textup{d}x\right)\mathclose{}\textup{d}s=w_{d-1}^{\rad}\lVert f\rVert_{L^1}\lVert \hat{f}\rVert_{L^1},
    \end{align*} so \begin{align*}
        \frac1{\lVert f\rVert_{L^1}\lVert \hat{f}\rVert_{L^1}}\int_{\mathbb{R}^d}|f|^2\le w_{d-1}^{\rad}.
    \end{align*} Taking the supremum over non-zero $f \in \mathcal{H}_{\rad}^d$ implies $w_d^{\rad}\le w_{d-1}^{\rad}$, as desired.
\end{proof}
Next, we introduce the restricted problem \begin{align}
    w_{0,d}^{\rad}\coloneqq \sup_{\substack{f \in \mathcal{H}_{\rad}^d\setminus \{\mathbf{0}\} \\ f(0)=\hat{f}(0)=0}}\frac{\lVert f\rVert_{L^2}^2}{\lVert f\rVert_{L^1}\lVert \hat{f}\rVert_{L^1}},
    \label{w_d0def}
\end{align} and show that \begin{lemma}
    With $w_{0,d}^{\rad}$ as defined above, we have $w_d^{\rad}\le 4w_{0,d}^{\rad}$ for all $d\ge 1$.
    \label{boundgenby0}
\end{lemma}
\begin{proof}
Let $f \in \mathcal{H}_{\rad}^d$.
    For each $\eps>0$, let $G_\eps(x)=G(x/\eps)$, where $G \in \mathcal{S}(\mathbb{R}^d)$ is defined by \begin{align*}
        \hat{G}(\xi)=\frac{2\pi |\xi|^2}{d}e^{-\pi |\xi|^2}.
    \end{align*} Note that $\hat{G}(0)=0$, and  \begin{align*}
    G(x)=-\frac1{2\pi d}\Delta(e^{-\pi |\cdot|^2})(x)=\frac1{2\pi d}(2\pi d-4\pi^2 |x|^2)e^{-\pi |x|^2},
\end{align*} so $G(0)=1$. Moreover, $\lVert \hat{G}\rVert_{L^1}=1$ since $\widehat{G}\ge 0$. Thus $G_\eps(0)=1$, $\hat{G_\eps}(0)=0$. Defining $\mathcal{E}_\eps$ by \begin{align}
    \mathcal{E}_\eps(x)=f(0)G_\eps(x)+\widehat{f}(0)\widehat{G_\eps}(x),
    \label{Epsxdef}
\end{align} we have \begin{align*}
    \mathcal{E_\eps}(0)&=f(0)G_\eps(0)+\widehat{f}(0)\widehat{G_\eps}(0)=f(0), \\
    \widehat{\mathcal{E_\eps}}(0)&=f(0)\widehat{G_\eps}(0)+\widehat{f}(0)G_\eps(0)=\widehat{f}(0).
\end{align*} Thus \begin{align}
    \left\lVert f-\mathcal{E}_\eps\right\rVert_{L^2}^2\le w_{0,d}^{\rad}\lVert f-\mathcal{E}_\eps\rVert_{L^1}\lVert \hat{f}-\hat{\mathcal{E}_\eps}\rVert_{L^1} \le w_{0,d}^{\rad}\mathopen{}\left(\lVert f\rVert_{L^1}+\lVert \mathcal{E}_\eps\rVert_{L^1}\right)\mathclose{}\mathopen{}\left(\lVert \hat{f}\rVert_{L^1}+\lVert \hat{\mathcal{E}_\eps}\rVert_{L^1}\right)\mathclose{}.
    \label{loveme}
\end{align} We now let $\eps \to 0$. For that, we note
\begin{align*}
    \left\lVert G_\eps\right\rVert_{L^2}^2=\int_{\mathbb{R}^d}\eps^d\left|G(x/\eps)\right|^2\textup{d}(x/\eps)=\eps^d\lVert G\rVert_{L^2}^2,
\end{align*} $\left\lVert G_\eps\right\rVert_{L^1}=\eps^d\lVert G\rVert_{L^1}$ and $ \left\lVert \widehat{G_\eps}\right\rVert_{L^1}=\lVert \hat{G}\rVert_{L^1}=1$ so  \begin{align}
\begin{split}
      \left\lVert \mathcal{E}_\eps\right\rVert_{L^1}&\le \eps^d|f(0)|\lVert G\rVert_{L^1}+|\hat{f}(0)|\le \eps^d|f(0)|\lVert G\rVert_{L^1}+\lVert f\rVert_{L^1} \\
     \left\lVert \widehat{\mathcal{E}_\eps}\right\rVert_{L^1}&\le \eps^d|\hat{f}(0)|\lVert G\rVert_{L^1}+|f(0)|\le \eps^d|\hat{f}(0)|\lVert G\rVert_{L^1}+\lVert \hat{f}\rVert_{L^1}.
\end{split}
  \label{partie}
\end{align} Letting $\eps \to 0$ in \eqref{loveme} gives, \begin{align*}
    \lVert f\rVert_{L^2}^2&\le 4w_{0,d}^{\rad}\lVert f\rVert_{L^1}\lVert \hat{f}\rVert_{L^1},
\end{align*} and the supremum over $f \in \mathcal{H}^d_{\rad}$ proves the lemma.
\end{proof}

Having established these Lemmas, it thus suffices to consider the restricted problem (where $f(0)=\hat{f}(0)=0$) in even dimensions. We also note we can further restrict to Schwartz functions, by using regular approximation combined with the subtraction trick employed above.

\section{Proof of the upper bound}\label{three}
As noted earlier, we take $\lambda \coloneqq \frac{d}{2} \in \mathbb{N}$. Let $f \in \mathcal{S}_{\rad}(\mathbb{R}^d)$, with $f(0)=\hat{f}(0)=0$. We abuse notation and let $f(r)$, $r >0$, denote the common value of $f(x)$ when $|x|=r$.
Following their notation, we define a holomorphic function $Z_f$, for such an $f$, on the strip $S_\lambda\coloneqq \left\{z \in \mathbb{C}: |\Im(z)|<\lambda\right\}$ by \begin{align*}
    Z_f(z)=\int_0^\infty r^{d-1}f(r)r^{-\lambda-iz}\textup{d}r=s_d^{-1}\int_{\mathbb{R}^d}f(x)|x|^{-\lambda-iz}\textup{d}x,
\end{align*} where \begin{align*}
    s_d=\frac{d\pi^{\frac{d}{2}}}{\Gamma(d/2+1)}
\end{align*} is the measure of the unit sphere in $\mathbb{R}^d$. We use the distributional Fourier transform formula \begin{align*}
    \hat{|\cdot|^{-\alpha}}(\xi)=\pi^{\alpha-\frac{d}{2}}\frac{\Gamma((d-\alpha)/2)}{\Gamma(\alpha/2)}|\xi|^{\alpha-d}, \quad \Re(\alpha) \in (0,d),
\end{align*} from \cite{steinreiz}, Chapter V, Section 1, Riesz potentials. So by symmetry, \begin{align*}
    Z_{f}(z)=s_d^{-1}\int_{\mathbb{R}^d}\hat{f}(\xi)\cdot \pi^{iz}\frac{\Gamma((\lambda-iz)/2)}{\Gamma((\lambda+iz)/2)}|\xi|^{-\lambda+iz}\textup{d}\xi= \pi^{iz} \frac{\Gamma((\lambda-iz)/2)}{\Gamma((\lambda+iz)/2)}Z_{\hat{f}}(-z),
\end{align*} so  \begin{align}
     Z_{f}(z)=\pi^{iz}\frac{\Gamma((\lambda-iz)/2)}{\Gamma((\lambda+iz)/2)}\cdot Z_{\hat{f}}(-z) 
     \label{eqnref}
\end{align} for all $z \in S_\lambda$. That $f(0)=\hat{f}(0)=0$ allows us to extend to $\Im(z)=\pm \lambda$.
Writing the integral for $Z_f(s)$ for real $s$ as \begin{align*}
    Z_f(s)=\int_0^\infty f(r)r^{\frac{d}{2}-1-is}\textup{d}r=\int_{\mathbb{R}} f(e^t)e^{\frac{td}{2}}e^{-its}\textup{d}t,
\end{align*} by Plancherel, we get \begin{align*}
    \int_{\mathbb{R}}|Z_f(s)|^2\textup{d}s=2\pi \int_{\mathbb{R}}|f(e^t)|^2e^{td}\textup{d}t=2\pi \int_{0}^\infty |f(r)|^2r^{d-1}\textup{d}r,
\end{align*} so \begin{align}
\lVert f\rVert_{L^2(\mathbb{R}^d)}^2=\frac{s_d}{2\pi}\int_{\mathbb{R}}|Z_f(s)|^2\textup{d}s.
    \label{goodeqn}
\end{align} If $z=s+i\lambda$, so $\lambda+iz=\lambda+is-\lambda=is$, then \begin{align*}
    \left|Z_f(s+i\lambda)\right|&=\left|\int_0^\infty r^{d-1}f(r)r^{-is} \textup{d}r\right|\le \frac{\lVert f\rVert_{L^1}}{s_d}.
\end{align*} For $z=s-i\lambda$, so $\lambda+iz=2\lambda+is$, $s\neq 0$, the reflection formula \eqref{eqnref} gives \begin{align*}
    \left|Z_{f}(s-i\lambda)\right|\le \pi^\lambda \left|\frac{\Gamma(-is/2)}{\Gamma(\lambda+is/2)}\right|\frac{\lVert \hat{f}\rVert_{L^1}}{s_d}\le \pi^\lambda \left|\frac{\Gamma(is/2)}{\Gamma(\lambda+is/2)}\right|\frac{\lVert \hat{f}\rVert_{L^1}}{s_d},
\end{align*} so (recall we are taking $\lambda \in \mathbb{N}$) \begin{align*}
    \log \left|Z_{f}(s-i\lambda)\right|&\le \lambda\log\pi+\log(\lVert \hat{f}\rVert_{L^1}/s_d)-\sum_{k=0}^{\lambda-1}\log\mathopen{}\left(\sqrt{k^2+s^2/4}\right)\mathclose{}.
\end{align*} 
Define $b_\lambda\colon \mathbb{R} \to \mathbb{R}\cup \{-\infty\}$ by \begin{align}
    b_\lambda(s)\coloneqq \sum_{k=0}^{\lambda-1}\log\mathopen{}\left(\sqrt{k^2+s^2/4}\right)\mathclose{}, \quad b_\lambda(0)=-\infty,
    \label{professional2}
\end{align} so the two bounds for the traces are \begin{align}
    \begin{split}
        \log\left|Z_f(s+i\lambda)\right|&\le \log(\lVert f\rVert_{L^1}/s_d), \\
    \log\left|Z_f(s-i\lambda)\right|&\le \lambda \log\pi +\log(\lVert \hat{f}\rVert_{L^1}/s_d)-b_\lambda(s).
    \end{split}
\label{models}
\end{align}
At this point, we majorise $\log|Z_f|$ by the harmonic extension of its traces; Theorem B.1 in \cite{CarvalhoCorso2026} implies one can majorise a bounded subharmonic function $u$ on the strip $\Im (z) \in (0,1)$, with continuous boundary traces $h=u(\cdot, 0)$, $g=u(\cdot, 1)$ by \begin{align*}
    U(x,y)&=(P_y*h)(x)+(P_{1-y}*g)(x), \quad
    P_y(x)=\frac1{2}\cdot \frac{\sin(\pi y)}{\cosh(\pi x)-\cos(\pi y)},
\end{align*}  assuming $h/\cosh(\pi\cdot), g/\cosh(\pi\cdot) \in L^1(\mathbb{R})$. In our case since $\log|Z_f|$ is bounded above, we can apply this principle to $\max(\log |Z_f|,-N)$, and let $N \to \infty$, applying monotone convergence to achieve our desired result.

For $u$ on the strip $S_\lambda$, one then has \begin{align*}
    u(2\lambda z-i\lambda)\le (P_y*u(2\lambda\cdot, -\lambda))(x)+(P_{1-y}*u(2\lambda \cdot, \lambda))(x).
\end{align*} at the midpoint $y=\frac1{2}$, so $ P_{\frac1{2}}(x)= \frac1{2\cosh(\pi x)}$,
we get \begin{align*}
    u(2\lambda x, 0)&\le \int_{\mathbb{R}}\frac1{2}\cdot \frac{u(2\lambda t,-\lambda)+u(2\lambda t,\lambda)}{\cosh(\pi (x-t))}\textup{d}t =\frac1{4\lambda }\int_{\mathbb{R}}\frac{u(s,-\lambda)+u(s,\lambda)}{\cosh\left(\pi \mathopen{}\left(\frac{2\lambda x}{2\lambda}-\frac{s}{2\lambda}\right)\mathclose{}\right)})\textup{d}s,
\end{align*} so \begin{align*}
      u(x,0)&\le \frac1{4\lambda }\int_{\mathbb{R}}\frac{u(s,-\lambda)+u(s,\lambda)}{\cosh\left(\pi \mathopen{}\left(\frac{x-s}{2\lambda}\right)\mathclose{}\right)}\textup{d}s=\frac1{2}\int_{\mathbb{R}}\frac{u(x-2\lambda t,-\lambda)+u(x-2\lambda t,\lambda)}{\cosh\left(\pi t\mathclose{}\right)}\textup{d}t.
\end{align*}
Applying this to $\log|Z_f|$ on the real line, we get
\begin{align}
    \log\left|Z_f(s)\right|&\le  \frac{\log(\lVert f\rVert_{L^1}/s_d)}{2}+\frac1{2}\int_{\mathbb{R}}\frac{\log|Z_f(s-2\lambda t-i\lambda)|}{\cosh(\pi t)}\textup{d}t.
    \label{kissland}
\end{align}
 To apply the upper bound in \eqref{models}, we need to verify $b_\lambda$ (shifted appropriately) is integrable against $1/\cosh(\pi t)$. From \eqref{professional2}, we see that the logarithmic singularity at $s=0$ is integrable, and $b_\lambda$ grows logarithmically in $s$. Thus, we may substitute the upper bound \eqref{models} in \eqref{kissland} to obtain
 \begin{align*}
    \log\left|Z_f(s)\right|&\le \frac{\lambda\log \pi+\log(\lVert f\rVert_{L^1}/s_d)+\log(\lVert \hat{f}\rVert_{L^1}/s_d)}{2}-\frac1{2}\int_{\mathbb{R}}\frac{b_\lambda(s-2\lambda t)}{\cosh(\pi t)}\textup{d}t,
\end{align*} thus implying \begin{align*}
    \left|Z_f(s)\right|^2\le \pi^\lambda s_d^{-2} \lVert f\rVert_{L^1}\lVert \hat{f}\rVert_{L^1}\exp\mathopen{}\left(-\int_{\mathbb{R}}\frac{b_\lambda(s-2\lambda t)}{\cosh(\pi t)}\textup{d}t\right)\mathclose{},
\end{align*} which we write as
\begin{align}
\begin{split}
    2\lambda s_d\left|Z_f(2\lambda s)\right|^2&\le 2\lambda \pi^\lambda s_d^{-1}\lVert f\rVert_{L^1}\lVert \hat{f}\rVert_{L^1}\exp(H_\lambda (s)),  \\
H_\lambda(s)&\coloneqq -\int_{\mathbb{R}}\frac{b_\lambda(2\lambda(s- t))}{\cosh(\pi t)}\textup{d}t.
\end{split}
 \label{mainineq3}
\end{align} 
 From \eqref{professional2}, \begin{align*}
   b_\lambda (2\lambda t)&=\sum_{k=0}^{\lambda-1}\log(\sqrt{k^2+\lambda^2t^2})=\lambda\log\lambda+\sum_{k=0}^{\lambda-1}\log(\sqrt{(k/\lambda)^2+t^2})\\
   &=\lambda\log\lambda+\sum_{k=1}^{\lambda}\log(\sqrt{(k/\lambda)^2+t^2})+\frac1{2}\log( t^2)-\frac1{2}\log(t^2+1),
\end{align*} and since $\log\sqrt{x^2+t^2}$ is increasing in $x\ge 0$, we get \begin{align*}
    b_{\lambda}(2\lambda t)\ge \lambda \log \lambda-\frac1{2}\log\mathopen{}\left(1+\frac1{t^2}\right)\mathclose{}+\lambda \int_{0}^{1}\log(\sqrt{x^2+t^2})\textup{d}x,
\end{align*} thus \begin{align*}
    -\int_{\mathbb{R}}\frac{b_\lambda(2\lambda(s-t))}{\cosh(\pi t)}\textup{d}t&\le -\lambda\log\lambda+\frac1{2}\int_{\mathbb{R}}\frac{\log(1+1/(t-s)^2)}{\cosh(\pi t)}\textup{d}t\\
    &-\lambda\int_{\mathbb{R}}\int_0^1\frac{\log(\sqrt{x^2+(t-s)^2})}{\cosh(\pi t)}\textup{d}x\textup{d}t.
\end{align*} As the convolution of two radially decreasing functions is also radially decreasing, we get \begin{align}
H_\lambda(s)&\le  -\lambda\log \lambda+E-\lambda\int_{\mathbb{R}}\int_0^1\frac{\log(\sqrt{x^2+t^2})}{\cosh(\pi t)}\textup{d}x\textup{d}t\eqcolon -\lambda\log \lambda+E+\lambda J_0 ,\label{ineq0} \\
    H_\lambda(s)&\le -\lambda\log \lambda+E-\lambda\int_{\mathbb{R}}\frac1{\cosh(\pi t)}\int_0^1\log(\sqrt{x^2+(t-s)^2})\textup{d}x\textup{d}t.
    \label{ineq}
\end{align} So \begin{align}
J_0=-\frac1{2}\int_{\mathbb{R}}\int_0^1\frac{\log(x^2+t^2)}{\cosh(\pi t)}\textup{d}x\textup{d}t.
    \label{J0def}
\end{align}
Finally, we need decay bounds on the convolution. We note that \begin{align*}
    \int_{\mathbb{R}}\frac{\log|s-t|}{\cosh(\pi t)}\textup{d}t &\ge \int_{|t|\le 1}\frac{\log|s-t|}{\cosh(\pi t)}\textup{d}t+\int_{|t|\le 1}\frac{\log |t|}{\cosh(\pi(s-t))}\textup{d}t.
\end{align*} If $|s|$ is large enough, this gives \begin{align*}
      \int_{\mathbb{R}}\frac{\log|s-t|}{\cosh(\pi t)}\textup{d}t&\ge C\log(|s|/2)-\frac{2}{\cosh(\pi |s|/2)}\ge C'\log |s|,
\end{align*} for constants $C,C' \in (0,1)$. So if $|s|$ is big enough, we have for some $c \in (0,1)$,  \begin{align*}
    \int_{\mathbb{R}}\frac1{\cosh(\pi t)}\int_0^1\log(\sqrt{x^2+(s-t)^2})\textup{d}x\textup{d}t\ge c\log |s|,
\end{align*} which implies \begin{align*}
        H_\lambda(s)\le -\lambda\log\lambda+E-c\lambda \log|s|,
    \end{align*} for all $s$ large enough (independent of $\lambda$).
Combining this with the universal bound \eqref{ineq0}, which says $ \exp(H_\lambda (s))\le e^E \lambda^{-\lambda }e^{\lambda J_0} $, gives \begin{align*}
   \int_{\mathbb{R}}\exp(H_\lambda(s))\textup{d}s&=\int_{|s|\le B}\exp(H_\lambda(s))\textup{d}s+\int_{|s|\ge B}\exp(H_\lambda(s))\textup{d}s \\
   &\le 2Be^E\lambda^{-\lambda } e^{\lambda J_0} +2\lambda^{-\lambda} e^E \int_{B}^\infty |s|^{-c\lambda}\textup{d}s, 
\end{align*} for large enough $B$. Letting $B$ (taking $\lambda$ big enough so $c\lambda>1$) be large enough gives 
 \begin{align*}
        \int_{\mathbb{R}}\exp(H_\lambda(s))\textup{d}s&\lesssim e^E\lambda^{-\lambda}e^{\lambda J_0}.
    \end{align*} 
Plugging this into \eqref{mainineq3} gives \begin{align*}
    s_d\int_{\mathbb{R}}|Z_f(s)|^2\textup{d}s\lesssim \lambda \lVert f\rVert_{L^1}\lVert \hat{f}\rVert_{L^1}\pi^\lambda s_d^{-1}\lambda^{-\lambda }e^{\lambda J_0},
\end{align*} so from \eqref{goodeqn}, \begin{align*}
    \frac1{\lVert f\rVert_{L^1}\lVert \hat{f}\rVert_{L^1}}\int_{\mathbb{R}^d}|f|^2\lesssim \lambda s_d^{-1}\pi^\lambda\lambda^{-\lambda}e^{\lambda J_0}.
\end{align*}
 From Stirling's approximation, we have \begin{align*}
    s_d^{-1}=\frac1{d}\frac{\Gamma(\lambda+1)}{\pi^\lambda}\sim \frac{\sqrt{2\pi \lambda }}{d}\cdot \pi^{-\lambda }(\lambda/e)^\lambda=\sqrt{\frac{\pi}{2\lambda}}\pi^{-\lambda}(\lambda/e)^\lambda,
\end{align*} so we get \begin{align*}
      \frac1{\lVert f\rVert_{L^1}\lVert \hat{f}\rVert_{L^1}}\int_{\mathbb{R}^d}|f|^2\lesssim \sqrt{\lambda} e^{\lambda(J_0-1)}.
\end{align*} This holds for any $f \in \mathcal{S}_{\rad}(\mathbb{R}^d)$ with $f(0)=\hat{f}(0)=0$. By density, we get for even $d$ large enough, 
    $w_{0,d}^{\rad}\lesssim \sqrt{d}e^{\frac{d(J_0-1)}{2}}$. The smaller $d$ cases follow from the fact that $w_d^{\rad}\le 1$. Combining this bound with Lemmas \ref{firstlemma} and \ref{boundgenby0}, we get \begin{align*}
    w_d^{\rad} \lesssim \sqrt{d}e^{\frac{d(J_0-1)}{2}},
\end{align*} for all dimensions $d$. To complete the proof of the upper bound, it remains to show that \begin{align}
    \frac{J_0-1}{2}=\frac{\log 2}{2}-\frac{2G}{\pi}.
    \label{qual}
\end{align}
\begin{proof}[Proof of \eqref{qual}]
Let \begin{align*}
    I(x)\coloneqq \frac1{2}\int_{\mathbb{R}}\frac{\log(x^2+t^2)}{\cosh(\pi t)}\textup{d}t
\end{align*} for $x>0$, so from \eqref{J0def} \begin{align}
J_0=-\int_0^1I(x)\textup{d}x.
    \label{J_0char}
\end{align} Observe that \begin{align*}
    I'(x)=\int_{\mathbb{R}}\frac{x}{x^2+t^2}\sech(\pi t)
\textup{d}t=\int_{\mathbb{R}}\frac1{x}\cdot \frac1{(t/x)^2+1}\sech(\pi t)\textup{d}t,
\end{align*} and using the fact that $\sech(\pi t)$ is a Fourier $1$-eigenfunction, we get \begin{align*}
    I'(x)=\int_{\mathbb{R}} \frac{e^{-2\pi x|\xi|}}{\cosh(\pi \xi)}\pi \textup{d}\xi=2\cdot 2\int_0^\infty \frac{e^{-2\pi x \xi}}{e^{\pi \xi}+e^{-\pi \xi}}\pi \textup{d}\xi=2\int_0^\infty \frac{e^{-2\pi  \xi(x+1/2)}}{1+e^{-2\pi \xi}}e^{2\pi \xi }2\pi e^{-2\pi\xi} \textup{d}\xi,
\end{align*} so \begin{align*}
    I'(x)=2\int_0^1\frac{t^{x+\frac1{2}-1}}{t+1}\textup{d}t\eqcolon 2\beta(x+1/2).
\end{align*} From \cite[entries 8.370,8.371, p.~907]{tableintergrals7th}, we have \begin{align*}
    \beta(x)=\frac1{2}\mathopen{}\left(\psi\mathopen{}\left(\frac{x+1}{2}\right)\mathclose{}-\psi\mathopen{}\left(\frac{x}{2}\right)\mathclose{}\right)\mathclose{},
\end{align*} where $\psi$ is the digamma function, so \begin{align*}
    I(x)=2\log(\Gamma)\mathopen{}\left(\frac{x}{2}+\frac3{4}\right)\mathclose{}-2\log(\Gamma)\mathopen{}\left(\frac{x}{2}+\frac1{4}\right)\mathclose{}+C,
\end{align*}  for some constant $C$. Note \begin{align*}
    I(0)=\int_{\mathbb{R}}\frac{\log(|t|)}{\cosh(\pi t)}\textup{d}t=\frac{2}{\pi}\int_0^\infty\frac{\log(\pi t)-\log(\pi)}{\cosh(\pi t)}\pi\textup{d}t=\frac2{\pi}\int_0^\infty\frac{\log(t)}{\cosh(t)}\textup{d}t-\log(\pi).
\end{align*}
From \cite[entry 4.371.1, p.~578]{tableintergrals7th}, we have \begin{align*}
 \int_0^\infty\frac{\log (t)}{\cosh( t)}\textup{d}t=\frac{\pi\log(2\pi)}{2}+\pi\mathopen{}\left(\log\Gamma(3/4)-\log\Gamma(1/4)\right)\mathclose{},
\end{align*} so we may deduce \begin{align*}
    I(x)=2\log(\Gamma)\mathopen{}\left(\frac{x}{2}+\frac3{4}\right)\mathclose{}-2\log(\Gamma)\mathopen{}\left(\frac{x}{2}+\frac1{4}\right)\mathclose{}+\log(2).
\end{align*} Integrating from $0$ to $1$, we get \begin{align*}
    \int_0^1I(x)\textup{d}x&=4\int_{\frac{3}{4}}^{\frac5{4}}\log\Gamma(x)\textup{d}x-4\int_{\frac1{4}}^{\frac3{4}}\log\Gamma(x)\textup{d}x+\log(2) \\
    &=4\int_{\frac{3}{4}}^{1}\log\Gamma(x)\textup{d}x+4\int_{0}^{\frac1{4}}\log\Gamma(1+x)\textup{d}x-4\int_{\frac1{4}}^{\frac3{4}}\log\Gamma(x)\textup{d}x+\log(2),
\end{align*} so from the Gamma recurrence, we get \begin{align*}
    \int_0^1I(x)\textup{d}x&=-\log(2)-1+4\int_{0}^{\frac1{4}}\log\Gamma(x)\textup{d}x+4\int_{\frac3{4}}^{1}\log\Gamma(x)\textup{d}x-4\int_{\frac1{4}}^{\frac3{4}}\log\Gamma(x)\textup{d}x \\
    &=-\log(2)-1+4\int_0^1\log\Gamma(x)\textup{d}x-8\int_{\frac1{4}}^{\frac3{4}}\log\Gamma(x)\textup{d}x.
\end{align*} Applying the reflection formula gives \begin{align*}
    \int_0^1I(x)\textup{d}x&=-\log(2)-1+2\int_{0}^{1}\log(\pi/\sin(\pi t))\textup{d}t-4\int_{\frac1{4}}^{\frac3{4}}\log(\pi/ \sin(\pi t))\textup{d}t \\
    &=-\log(2)-1+4\int_{\frac1{4}}^{\frac3{4}}\log(\sin(\pi t))\textup{d}t-2\int_0^1\log(\sin(\pi t))\textup{d}t,
\end{align*} which we write as \begin{align*}
    \int_0^1I(x)\textup{d}x=-\log(2)-1+2\int_0^1\log(\sin(\pi t))\textup{d}t-8\int_0^{\frac1{4}}\log(\sin(\pi t))\textup{d}t.
\end{align*}
For the first integral, we have
\begin{align*}
2\int_0^1\log(\sin(\pi t))\textup{d}t=\int_0^1\log\mathopen{}\left(\frac{2-2\cos(2\pi t)}{4}\right)\mathclose{}\textup{d}t=-2\log 2,
\end{align*} so \begin{align*}
    \int_0^1I(x)\textup{d}x=-3\log(2)-1-8\int_{0}^{\frac1{4}}\log(\sin(\pi t))\textup{d}t=-3\log(2)-1-\frac{8}{\pi}\int_{0}^{\frac{\pi}{4}}\log(\sin(t))\textup{d}t.
\end{align*} From \cite[entry 4.224.2, p.~530]{tableintergrals7th},  we have \begin{align}
    \int_{0}^{\frac{\pi}{4}}\log(\sin(t))\textup{d}t=-\frac{\pi}{4}\log2-\frac{G}{2},
    \label{catalanchar}
\end{align}
 thus \begin{align*}
    \int_0^1I(x)\textup{d}x=-\log2-1+\frac{4G}{\pi}.
\end{align*} From \eqref{J_0char}, this gives $J_0=\log 2+1-4G/\pi$, as desired.

\end{proof}

\section{Proof of the lower bound }\label{lowerbound}
Given a real-valued integrable radial function $h$ on $\mathbb{R}^d$, define $H_h\colon \overline{\mathbb{D}}\setminus \{-1\}\to \mathbb{C}$ by  \begin{align*}
    H_h(z)=(1-z)^{\frac{d}{2}}\int_{\mathbb{R}^d}h(x)e^{-\frac{1-z}{1+z}\pi |x|^2}\textup{d}x.
\end{align*} Observe the map $z \mapsto \frac{1-z}{1+z}$ maps the disk to the right-half plane: \begin{align*}
    2\Re\mathopen{}\left(\frac{1-z}{1+z}\right)\mathclose{}=\frac{1-z}{1+z}+\frac{1-\overline{z}}{1+\overline{z}}=\frac{1+\overline{z}-z-|z|^2+1-\overline{z}+z-|z|^2}{|1+z|^2}=\frac{2-2|z|^2}{|1+z|^2}\ge 0,
\end{align*} if $|z|\le 1$. By symmetry we have for $z \in \mathbb{D}$, \begin{align*}
    H_{\hat{h}}(z)&=(1-z)^{\frac{d}{2}}\int_{\mathbb{R}^d}\hat{h}(x)e^{-\frac{1-z}{1+z}\pi |x|^2}\textup{d}x=(1-z)^{\frac{d}{2}}\mathopen{}\left(\frac{1+z}{1-z}\right)^{\frac{d}{2}}\int_{\mathbb{R}^d}h(x)e^{-\frac{1+z}{1-z}\pi |x|^2}\textup{d}x =H_h(-z),
\end{align*} thus $H_h$ is even if and only if $h$ is a $1$-eigenfunction; the reverse implication follows from taking $z \in \mathbb{D}$ and using the injectivity of the Laplace transform.

We will now construct an even holomorphic function $H$ on $\mathbb{D}$ extending continuously to $\partial \mathbb{D}$, for which there exists a 1-eigenfunction $h \in \mathcal{H}_{\rad}^d$ such that $H_h=H$. This will be our candidate for lower bounding $w_d^{\rad}$.

Let $\mathcal{H}$ be the operator acting on real-valued square-integrable functions $v\colon \T \to \mathbb{R}$, defined by \begin{align*}
    \mathcal{H}v(z)=\hat{v}(0)+2\sum_{n=1}^\infty \hat{v}(n)z^n.
\end{align*} So $\mathcal{H}$ gives a holomorphic function $V$ on $\mathbb{D}$ with boundary trace $\Re(V)=v$ in $L^2(\T$).
 Let \begin{align}
    L=\mathcal{H}(f), \qquad f(z)=\log\min(|1-z|,|1+z|) \quad \forall z \in \T,
    \label{Ldef}
\end{align} (with $-\infty$ at $\pm 1$) so $L$ is holomorphic on $\mathbb{D}$. We will show that \begin{align}
L(z)=\hat{f}(0)+\frac{\log(1-z^2)}{2}-\frac4{\pi}\int_0^z\frac{\arctan(w)}{1-w^2}\textup{d}w, 
    \label{lform}
\end{align} with the principal branches. In particular $L$ extends to a holomorphic function on $\mathbb{C}\setminus \mathopen{}\left(\pm [1,\infty)\cup \pm i[1,\infty)\right)\mathclose{}$.
Thus $L$ extends continuously to $\partial \mathbb{D} \simeq \T$ away from $\pm 1$, $\pm i$ , with \begin{align*}
    \Re(L)(z)=\log \min(|1-z|,|1+z|), \quad z \in \mathbb{T}\setminus \{\pm 1\}.
\end{align*} So $\Re(L)$ is the Poisson extension of $\log \min(|1-z|,|1+z|)$. Also note $L$ is even since $f$ is even, so its Fourier frequencies consist of the even ones.
Then the function $H\colon \mathbb{D} \to \mathbb{C}$ defined by \begin{align}
    H(z)\coloneqq \exp\mathopen{}\left(dL(z)/2\right)\mathclose{}
    \label{myname}
\end{align} is even, and holomorphic on the disk. Applying the mean-value theorem to $\log |H|$ and taking the limit as $r \to 1^{-}$, we have \begin{align*}
    \begin{split}
        \log |H(0)|&=\frac1{2\pi}\int_0^{2\pi}\log|H(e^{ i \theta})|\textup{d}\theta=\frac1{2\pi}\cdot \frac{d}{2}\int_{-\pi}^{\pi}\log\min(|1-e^{ i \theta}|,|1+e^{ i \theta}|)\textup{d}\theta \\
    &=\frac{1}{2\pi}\cdot \frac{d}{2}\cdot 2\int_{-\frac{\pi}{2}}^{\frac{\pi}{2}}\log|1-e^{i\theta}|\textup{d}\theta=\frac1{2\pi}\cdot \frac{d}{2}\int_{-\frac{\pi}{2}}^{\frac{\pi}{2}}\log(4\sin^2(\theta/2))\textup{d}\theta 
    \end{split}
\end{align*} so \begin{align*}
    \log|H(0)|=\frac{d}{\pi}\int_0^{\frac{\pi}{4}}\log(4\sin^2(t))\textup{d}t=\frac{2d}{\pi}\int_0^{\frac{\pi}{4}}\log(2\sin(t))\textup{d}t.
\end{align*} From \eqref{catalanchar}, this gives \begin{align}
    \log|H(0)|=-\frac{dG}{\pi}.
    \label{blacko}
\end{align}
We now bound the $L^1$-norm of the induced $h$. That is, we wish to establish the existence of an integrable function $h \in L^1_{\rad}(\mathbb{R}^d)$ such that $H_h=H$; that $h$ is a $1$-eigenfunction will follow from the fact that $H$ is even, as established earlier.
Note
\begin{align}
    \frac{1-e^{i\alpha}}{1+e^{i\alpha}}=\frac{(1-e^{i\alpha})(1+e^{-i\alpha})}{|1+e^{-i\alpha}|^2}=\frac{e^{-i\alpha}-e^{i\alpha}}{2+2\cos(\alpha)}=-\frac{2i\sin(\alpha)}{2+2\cos(\alpha)}=-i\tan(\alpha/2).
    \label{heart}
\end{align} So if $h$ satisfies $H_h=H$, $H$ as in \eqref{myname},  then we must have for all $\alpha \in (-\pi ,\pi)$,
 \begin{align}
    \exp\mathopen{}\left(\frac{dL(e^{i\alpha})}{2}\right)\mathclose{}=(1-e^{i\alpha})^{\frac{d}{2}}\int_{\mathbb{R}^d}h(x)e^{i\tan(\alpha/2)\pi |x|^2}\textup{d}x.
    \label{identity}
\end{align} This is sufficient; if \eqref{identity} holds, then $H_h-H$ is a holomorphic function on $\mathbb{D}$, with $0$ trace (at least away from $\pm 1,\pm i$), thus it must be identically $0$. We now turn the right-hand side of \eqref{identity} into a Fourier transform.
Converting to polar coordinates, \begin{align}
   \begin{split}
        &\int_{\mathbb{R}^d}h(x)e^{i\tan(\alpha/2)\pi |x|^2}\textup{d}x= \int_{0}^\infty s_dr^{d-1}h(r)e^{i\tan(\alpha/2)\pi r^2}\textup{d}r \\
    &=\frac1{2\pi}\int_{0}^\infty s_dr^{d-2}h(r)e^{i\tan(\alpha/2)\pi r^2}2\pi r\textup{d}r =\frac1{2\pi}\int_0^\infty s_d\mathopen{}(\sqrt{{t}/{\pi}})^{d-2}\mathclose{}h(\sqrt{t/\pi})e^{i\tan(\alpha/2)t}\textup{d}t.
   \end{split}
   \label{babe2}
\end{align} Let \begin{align}
    h_\mathbf{f}(t)=s_d\mathbbm{1}_{[0,\infty)}(t)(t_+/\pi)^{\frac{d}{2}-1}h(\sqrt{t_+/\pi}),
    \label{lonelystar}
\end{align}
so \eqref{identity} becomes \begin{align*}
    \exp\mathopen{}\left(\frac{dL(e^{i\alpha})}{2}\right)\mathclose{}= (1-e^{i\alpha})^{\frac{d}{2}}(h_\mathbf{f})^\vee(\tan(\alpha/2)), \quad \forall \alpha \in (-\pi,\pi)
\end{align*} which we write as \begin{align}
    h_\mathbf{f}^\vee(\xi)=U\mathopen{}\left(e^{2i\arctan(\xi)}\right)^{\frac{d}{2}}\mathclose{}\coloneqq \mathopen{}\left(1-e^{2i\arctan(\xi)}\right)^{-\frac{d}{2}}\mathclose{}\exp\mathopen{}\left(\frac{dL(e^{2i\arctan(\xi)})}{2}\right)\mathclose{},
    \label{129}
\end{align} so \begin{align}
    U(z)=\frac{\exp(L(z))}{1-z}.
    \label{hongkong}
\end{align} As before, let $\lambda=\frac{d}{2}$. We want to solve for an integrable function $h_\mathbf{f}$ supported on the right-half line, for \eqref{129} holds. Recall from \eqref{heart} \begin{align*}
    \frac{1-e^{2i\arctan(\xi)}}{1+e^{2i\arctan(\xi)}}=-i\xi\eqcolon \frac{1-z}{1+z},
\end{align*} and as \begin{align*}
    \frac{1-\frac{1-z}{1+z}}{1+\frac{1-z}{1+z}}=\frac{1+z-1+z}{1+z+1-z}=\frac{2z}{2}=z,
\end{align*} we get \begin{align}
   z=e^{2i\arctan(\xi)}=\frac{1+i\xi}{1-i\xi}\eqcolon m(\xi).
   \label{thetown}
\end{align} We then write \eqref{129} as \begin{align}
h_\mathbf{f}^\vee(\xi)=U(m(\xi))^{\frac{d}{2}}.
    \label{29}
\end{align}
 Now we establish that $h_\mathbf{f}$ solving \eqref{29} is supported on the right-half line. Write $m(\xi)$ as $m(\xi)=\frac2{1-i\xi}-1$, so  \begin{align}
    m'(\xi)=\frac{2i}{(1-i\xi)^2}=\frac{i(1+m(\xi))^2}{2}.
    \label{adaptation}
\end{align} For fixed $y>0$, we have  \begin{align*}
    \int_{\mathbb{R}}|U^{\lambda}(m(x+iy))|^2\textup{d}x&= \int_{\mathbb{R}}|U(m(x+iy))|^{2\lambda}m'(x+iy)\cdot \frac2{i(1+m(x+iy))^2}\textup{d}x,
\end{align*} so \begin{align}
 \int_{\mathbb{R}}|U^{\lambda}(m(x+iy))|^2\textup{d}x&=\int_{\gamma_y}\frac{2|U(z)|^{2\lambda}}{i(1+z)^2}\textup{d}z,
    \label{sine}
\end{align}
where $\gamma_y$ is the circle $\{m(x+iy): x \in \mathbb{R}\}$. Now, since $\Re(L)$ is the Harmonic extension of $\log\min(|1-z|,|1+z|)$, we have on $\T$, 
\begin{align*}
    \Re(L(z))\le \log|1-z|=\frac{\log|1-z|^2}{2}=\Re(\log(1-z)),
\end{align*} hence $|U(z)|\le 1$ for all $z \in \mathbb{D}$. We also have \begin{align*}
    \Re(L(z))&\le \min\mathopen{}\left(\Re(\log(1-z)),\Re(\log(1+z))\right)\mathclose{}\\
    &=\min(\log|1-z|,\log|1+z|)\le \log|1-z|+\log|1+z|,
\end{align*} since the triangle inequality implies one of $|1-z|,|1+z|$ is at least $1$. So we have \begin{align*}
    |U(z)|\le 1,|1+z| \quad \forall z \in \mathbb{D},
\end{align*} thus for all $\lambda\ge 1$,
 \begin{align}
    \begin{split}
       \left|\frac{U(z)^\lambda}{1+z}\right|\le 1.
    \end{split}
    \label{align}
\end{align} From \eqref{sine}, this implies
 \begin{align}
    \sup_{y>0} \int_{\mathbb{R}}|U^\lambda(m(x+iy))|^2\textup{d}x\le  2\sup_{y >0}|\gamma_y|=4\pi,
    \label{lushlife}
\end{align} so by Paley-Wiener, this implies $h_\mathbf{f}$ is square-integrable, and supported on the right-half line. 

Now, let $h$ be the radial function on $\mathbb{R}^d$ defined implicitly by \eqref{lonelystar}. Note from \eqref{babe2}, \begin{align}
    \left\lVert h\right\rVert_{L^1(\mathbb{R}^d)}=\frac{\left\lVert h_\mathbf{f}\right\rVert_{L^1(\mathbb{R})}}{2\pi}
    \label{hong}.
\end{align} Once we show that $h_\mathbf{f}$ is integrable, we will get that $h$ is integrable (that $h$ is a $1$-eigenfunction follows from evenness), thus \eqref{identity} would determine $H_h(z)$ for $z \in \partial \mathbb{D}\setminus \{ -1\}$. Since $h \in L^1(\mathbb{R}^d)$, both $H_h$ and $H$ have well-defined traces on $\partial \mathbb{D}\setminus \{\pm 1,\pm i\}$ (once we prove \eqref{lform}). Combining this with \eqref{29} gives $H_h=H$ on $\partial \mathbb{D}\setminus \{\pm 1,\pm i\}$. We also note that $H_h$ is bounded on $\mathbb{D}$, $|H_h(z)|\le 2^{\frac{d}{2}}\lVert h\rVert_{L^1(\mathbb{R}^d)}$, as is $H$ (the boundedness of $H$, defined in \eqref{myname}, follows from the fact that $\Re(L)$ is the Poisson extension of $\log\min(|1-z|,|1+z|)$, which is bounded above on $\T$). By considering the difference $H_h-H$ combined with uniqueness, we get $H_h=H$.

Once we have $H_h=H$, this would imply $h$ is real-valued. Once we show \eqref{lform}, this combined with \eqref{myname} would imply $H(\overline{z})=\overline{H(z)}$, thus $H_h(\overline{z})=\overline{H_h(z)}$. This implies $H_h=H_{\overline{h}}$ and combining with the injectivity of the Laplace transform, since $h$ is radial, implies $h$ is real-valued.\newline

We need to bound the $L^1$ norm of $h_\mathbf{f}$. Write the expression for $h_\mathbf{f}^\vee$ in \eqref{29} as \begin{align*}
    h_\mathbf{f}^\vee(\xi)=F(\xi)^{\lambda}, \quad F(\xi)\coloneqq U(m(\xi)).
\end{align*}  To estimate $\lVert h_\mathbf{f}\rVert_{L^1}$, we use Sobolev embedding, \begin{align*}
    \left\lVert h_\mathbf{f}\right\rVert_{L^1}^2\lesssim \lVert h_\mathbf{f}^\vee \rVert_{L^2}^2+\lVert (h_\mathbf{f}^\vee)'\rVert_{L^2}^2.
\end{align*} From \eqref{lushlife}, we have $\lVert h_\mathbf{f}^\vee\rVert_{L^2}^2\lesssim 1$ uniformly in $d$. Note \begin{align*}
    (h_\mathbf{f}^\vee)'(\xi)=\lambda F(\xi)^{\lambda-1}F'(\xi),
\end{align*} and since $|F(\xi)|\le 1$, we have \begin{align*}
    \lVert h_\mathbf{f}\rVert_{L^1}^2&\lesssim 1+\lambda^2\int_{\mathbb{R}}|F'(\xi)|^2\textup{d}\xi=1+\lambda^2\int_{\mathbb{R}}|U'(m(\xi))|^2|m'(\xi)|^2\textup{d}\xi \\
    &=1+\lambda^2\int_{\mathbb{R}}|U'(m(\xi))|^2\overline{m'(\xi)}m'(\xi)\textup{d}\xi
\end{align*} (Here the derivative is that of $U$ on the torus). From \eqref{adaptation}, we have \begin{align*}
    \overline{m'(\xi)}=-\frac{i}{2}\mathopen{}\left(1+\frac1{m(\xi)}\right)\mathclose{}^2,
\end{align*} since $|m(\xi)|=1$.
Writing the upper bound as  \begin{align*}
    \left\lVert h_\mathbf{f}\right\rVert_{L^1}^2&\lesssim 1-\frac{i\lambda^2}{2}\bigintsss_{\mathbb{R}}|U'(m(\xi))|^2(1+1/m(\xi))^2m'(\xi)\textup{d}\xi \\
    &=1-\frac{i\lambda^2}{2}\int_{\T}|U'(z)|^2(1+z^{-1})^2\textup{d}z,
\end{align*} gives \begin{align}
    \lVert h_\mathbf{f}\rVert_{L^1}^2\lesssim 1+2\lambda^2\int_{\T}|U'(z)|^2|\textup{d}z|.
    \label{floqo}
\end{align}
Recall from \eqref{hongkong}, \begin{align}
    U(z)=\frac{e^{L(z)}}{1-z}=\exp\mathopen{}\left(L(z)-\log(1-z)\right)\mathclose{}.
    \label{four}
\end{align}  We now calculate $L$, which is given by \begin{align*}
    L(z)=\mathcal{H}(f)(z)=\hat{f}(0)+2\sum_{n>0}\hat{f}(n)z^n.
\end{align*}
\begin{proof}[Proof of \eqref{lform}]
    Recall $f(z)=\log\min(|1-z|,|1+z|)$,  so $f$ has only even Fourier frequencies. 
To that end, for $n>0$ even, we have
\begin{align*}
    2\pi\hat{f}(n)&=\int_{-\frac{\pi}{2}}^{\frac{\pi}{2}}f(e^{i\theta})e^{-in\theta}\textup{d}\theta+\int_{\frac{\pi}{2}}^{\frac{3\pi}{2}}f(e^{i\theta})e^{-in\theta}\textup{d}\theta \\
    &=\int_{-\frac{\pi}{2}}^{\frac{\pi}{2}}f(e^{i\theta})e^{-in\theta}\textup{d}\theta+\int_{-\frac{\pi}{2}}^{\frac{\pi}{2}}f(e^{i(\pi+t)})e^{-in(\pi+t)}\textup{d}t=2\int_{-\frac{\pi}{2}}^{\frac{\pi}{2}}f(e^{i\theta})e^{-in\theta}\textup{d}\theta,
\end{align*} since $f$ is even. For $\theta \in [-\pi/2,\pi/2]$, $f\mathopen{}\left(e^{i\theta}\right)\mathclose{}=\log|1-e^{i\theta}|$, so  \begin{align*}
       -2\pi \hat{f}(n)&=\int_{-\frac{\pi}{2}}^{\frac{\pi}{2}}-\log|1-e^{i\theta}|^2e^{-in\theta}\textup{d}\theta.
\end{align*}
 Using the fact that, $\log|1-w|^2=\log(1-w)+\log(1-\overline{w})=2\Re\mathopen{}\left(\log(1-w)\right)\mathclose{}$ and expanding the latter in Taylor series, we get 
\begin{align*}
    -2\pi \hat{f}(n)&=\int_{-\frac{\pi}{2}}^{\frac{\pi}{2}}\sum_{k=1}^\infty \frac{2\cos(k\theta)e^{-in\theta}}{k}\textup{d}\theta=\int_{-\frac{\pi}{2}}^{\frac{\pi}{2}}\sum_{k=1}^\infty\frac{2\cos(k\theta)\cos(n\theta)}{k}\textup{d}\theta \\
    &=\sum_{k=1}^\infty\int_{-\frac{\pi}{2}}^{\frac{\pi}{2}} \frac{\cos((k-n)\theta)}{k}\textup{d}\theta+\sum_{k=1}^\infty\int_{-\frac{\pi}{2}}^{\frac{\pi}{2}} \frac{\cos((n+k)\theta)}{k}\textup{d}\theta.
\end{align*} Observe for $m\neq 0$ an even integer, one has $\int_{-{\pi}/{2}}^{{\pi}/{2}}\cos(m\theta)\textup{d}\theta=2\sin(\pi m/2)/m=0$. Thus \begin{align*}
    -2\pi \hat{f}(n)&=\frac{\pi }{n}+2\sum_{k>0\text{ odd }}\mathopen{}\left(\frac{(-1)^{\frac{k+n-1}{2}}}{k(n+k)}+\frac{(-1)^{\frac{k-n-1}{2}}}{k(k-n)}\right)\mathclose{}\\
    &=\frac{\pi }{n}+\frac2{n}\sum_{k>0\text{ odd }}\mathopen{}\left((-1)^{\frac{k+n-1}{2}}\mathopen{}\left(\frac1{k}-\frac1{n+k}\right)-(-1)^{\frac{k-n-1}{2}}\mathopen{}\left(\frac1{k}-\frac1{k-n}\right)\mathclose{}\right)\mathclose{},
\end{align*} and since $n$ is even, this reduces to \begin{align*}
    -2\pi \hat{f}(n)&=\frac{\pi}{n}+\frac2{n}\sum_{k>0 \text{ odd }}\mathopen{}\left(\frac{(-1)^{\frac{k-n-1}{2}}}{k-n}-\frac{(-1)^{\frac{k+n-1}{2}}}{k+n}\right)\mathclose{}=\frac{\pi}{n}+\frac2{n}\sum_{\substack{p \in [1-n, n] \\ p\text{ odd}}}\frac{(-1)^{\frac{p-1}{2}}}{p} \\
    &=\frac{\pi}{n}+\frac4{n}\sum_{\substack{p \in [1,n-1] \\ p\text{ odd}}}\frac{(-1)^{\frac{p-1}{2}}}p.
\end{align*} Recognising the sum as the partial sum of the Taylor series of $\arctan(z)$, we see that  \begin{align*}
    \sum_{\substack{p \in [1,n-1] \\ p\text{ odd}}}\frac{(-1)^{\frac{p-1}{2}}}p
\end{align*} is the coefficient of $z^{n-1}$ in the Taylor series expansion of $\frac{\arctan(z)}{1-z^2}$. Thus summing over even $n>0$, we get \begin{align*}
    2\sum_{n>0}\hat{f}(n)z^n=\frac{\log(1-z^2)}{2}-\frac4{\pi}\int_0^z\frac{\arctan(w)}{1-w^2}\textup{d}w,
\end{align*} as claimed.
\end{proof}

From \eqref{four}, we then have \begin{align*}
    U(z)&=\exp\mathopen{}\left(\hat{f}(0)+\frac{\log(1-z)+\log(1+z)}{2}-\log(1-z)-\frac4{\pi}\int_0^z\frac{\arctan(w)}{1-w^2}\textup{d}w\right)\mathclose{} \\
    &=\exp\mathopen{}\left(\hat{f}(0)+\frac{\log(1+z)-\log(1-z)}{2}-\frac4{\pi}\int_0^z\frac{\arctan(w)}{1-w^2}\textup{d}w\right)\mathclose{}, 
\end{align*} thus \begin{align*}
    U'(z)&=\mathopen{}\left(\frac1{2}\mathopen{}\left(\frac1{1+z}+\frac1{1-z}\right)\mathclose{}-\frac4{\pi}\cdot \frac{\arctan(z)}{1-z^2}\right)\mathclose{}U(z) =\mathopen{}\left(\frac1{1-z^2}-\frac4{\pi}\cdot \frac{\arctan(z)}{1-z^2}\right)\mathclose{}\frac{e^{L(z)}}{1-z} \\&=\frac4{\pi}\cdot \frac{\frac{\pi}{4}-\arctan(z)}{1-z}\cdot \frac{e^{L(z)}}{(1-z)(1+z)},
\end{align*} for $z \in \mathbb{D}$. 
For $z \in \T$, we have \begin{align*}
    \left|\frac{e^{L(z)}}{(1-z)(1+z)}\right|=\frac{\min(|1-z|,|1+z|)}{|(1-z)(1+z)|}=\max(|1-z|,|1+z|)^{-1}\le 1.
\end{align*} Moreover, since we defined $L$ in such a way that \begin{align*}
    \Re(L(z))-\Re(\log(1-z^2))
\end{align*} is the harmonic extension of its restriction to $\partial \mathbb{D}$, we may deduce
 $\left|e^{L(z)}/(1-z^2)\right|\le 1$ for all $z \in \mathbb{D}$. From the expression for $U'$, this implies \begin{align*}
    |U'(z)|\le \frac4{\pi}\cdot \frac{\left|\arctan(z)-\frac{\pi}{4}\right|}{|1-z|} \quad \forall z \in \mathbb{D}.
\end{align*} From the Taylor series \begin{align*}
    \arctan(z)=\sum_{m=0}^\infty \frac{(-1)^m}{2m+1}z^{2m+1},
\end{align*} we can see $\arctan(z)$ is square integrable on the torus. Combining this with the derivative  $\arctan'(z)=\frac1{1+z^2}$ which is bounded near $z=1$, we may deduce, from \eqref{hong} and \eqref{floqo}, \begin{align*}
   \left\lVert h\right\rVert_{L^1(\mathbb{R}^d)}^2= \frac{\left\lVert h_\mathbf{f}\right\rVert_{L^1(\mathbb{R})}^2}{4\pi^2}\lesssim d^2,
\end{align*} $d\ge 2$,
 as desired. So $\lVert h\rVert_{L^1}\le Cd$ for some constant $C$ (independent of $d$). From the definition of $H_h$, \begin{align*}
    \left|\int_{\mathbb{R}^d}h(x)e^{-\pi |x|^2}\textup{d}x\right|=|H_h(0)|=|H(0)|=\exp(\log|H(0)|)=\exp(-dG/\pi),
\end{align*} from \eqref{blacko}. Combining this with Cauchy-Schwarz gives  \begin{align*}
    \frac{\lVert h\rVert_{L^2}^2}{\lVert h\rVert_{L^1}^2}\gtrsim \frac1{d^2}\left|\int_{\mathbb{R}^d}h(x)\cdot 2^{\frac{d}{4}}e^{-\pi |x|^2}\textup{d}x\right|^2= \frac{2^{\frac{d}{2}}\exp(-2dG/\pi)}{d^2}=d^{-2}2^{\frac{d}{2}}e^{-\frac{2dG}{\pi}}\mathclose{},
\end{align*} for $d\ge 2$ initially. $d=1$ can be induced by testing against a Gaussian, decreasing the constant if necessary. This gives the lower bound in Theorem \ref{mainthm2}, as claimed.

 \section{Numerics}\label{dualform}
 Define $\k_d\colon \mathbb{R}^d\times \mathbb{R}^d \to \mathbb{R}$ by $\k_d(x,y)=\kappa_d(|x||y|)$. We note $\k_d$ is a $+1$-eigenfunction (distributionally). Indeed, for any $f,g \in \mathcal{S}(\mathbb{R}^d)$, one has \begin{align*}
    (f\otimes g, \k_d)&=\int_{\mathbb{R}^d}\int_{\mathbb{R}^d}f(x)g(y)\k_d(x,y)\textup{d}x\textup{d}y=\int_{\mathbb{R}^d}f(x)\mathopen{}\left(\int_{\mathbb{R}^d}\k_d(x,y)g(y)\textup{d}y\right)\mathclose{}\textup{d}x,
\end{align*} and using the fact that $\k_d(x,y)$ is radial in each of $x$ and $y$, we have \begin{align}
     (f\otimes g, \k_d)=\int_{\mathbb{R}^d}f(x)\mathopen{}\left(\int_{\mathbb{R}^d}\k_d(x,y)\mathcal{R}_dg(y)\textup{d}y\right)\mathclose{}\textup{d}x=\int_{\mathbb{R}^d}f(x)\mathcal{F}\mathcal{R}_dg(x)\textup{d}x,
     \label{wonderland}
\end{align} where $\mathcal{R}_d$ is the radialisation defined in \eqref{radialisationdef}. Using the fact that $\mathcal{F}$ commutes with $\mathcal{R}_d$ and the symmetry of the Fourier transform, we have \begin{align*}
     (f\otimes g, \k_d)=\int_{\mathbb{R}^d}\widehat{f}(x)\mathcal{R}_d g(x)\textup{d}x=\int_{\mathbb{R}^d}\widehat{f}(x)\mathcal{R}_d g(-x)\textup{d}x=\int_{\mathbb{R}^d}\widehat{f}(x)\mathcal{F}\mathcal{R}_d\hat {g}(x)\textup{d}x.
\end{align*} Comparing this to the right-hand side of \eqref{wonderland}, this implies $(f\otimes g,\k_d)=(\hat{f}\otimes \hat{g},\k_d)$.  This holds for all $f,g \in \mathcal{S}(\mathbb{R}^d)$, and since \begin{align*}
    \Span\mathopen{}\left\{f\otimes g: f,g \in \mathcal{S}(\mathbb{R}^d)\right\}
\end{align*} is dense in $\mathcal{S}\mathopen{}\left(\mathbb{R}^{2d}\right)\mathclose{}$, we obtain $\hat{\k_d}=\k_d$ distributionally, as claimed. 

We now note the following. For any $-1$-eigendistribution $u \in L^\infty(\mathbb{R}^{d}\times \mathbb{R}^d)$ that is also symmetric, so $u(x,y)=u(y,x)$ for all $x,y \in \mathbb{R}^d$,  we have for any $f \in \mathcal{S}_{\rad}(\mathbb{R}^d)$, \begin{align*}
    \int_{\mathbb{R}^d}\int_{\mathbb{R}^d}f(x)\hat{f}(y)u(x,y)\textup{d}x\textup{d}y&= \int_{\mathbb{R}^d}\int_{\mathbb{R}^d}\hat{f}(x)f(y)\hat{u}(x,y)\textup{d}x\textup{d}y =-\int_{\mathbb{R}^d}\int_{\mathbb{R}^d}\hat{f}(x)f(y)u(x,y)\textup{d}x\textup{d}y \\
    &=-\int_{\mathbb{R}^d}\int_{\mathbb{R}^d}\hat{f}(y)f(x)u(x,y)\textup{d}x\textup{d}y,
\end{align*} by symmetry, thus \begin{align*}
     \int_{\mathbb{R}^d}\int_{\mathbb{R}^d}f(x)\hat{f}(y)u(x,y)\textup{d}x\textup{d}y=0.
\end{align*}
Consequently,
\begin{align*}
    \frac{\lVert f\rVert_{L^2}^2}{\lVert f\rVert_{L^1}\lVert \hat{f}\rVert_{L^1}}=\frac1{\lVert f\otimes \hat{f}\rVert_{L^1}}\int_{\mathbb{R}^{2d}}(f\otimes \hat{f})(\k_d-u),
\end{align*} for any non-zero real-valued $f \in \mathcal{S}_{\rad}(\mathbb{R}^d)$, thus for any non-zero $f \in \mathcal{H}_{\rad}^d$, invoking the density of $\mathcal{S}_{\rad}(\mathbb{R}^d)$ in $\mathcal{H}_{\rad}^d$ under the $\mathcal{H}^d_{\rad}$ norm in \eqref{Hradnorm}. Applying H\"{o}lder's inequality, taking the supremum over $f \in \mathcal{H}_{\rad}^d\setminus \{\mathbf{0}\}$ gives \begin{align*}
    w_d^{\rad}\le \left\lVert \k_d-u\right\rVert_{L^\infty}.
\end{align*} How to construct such a function $u$? One way is as follows. Let $\gamma=\sum_{j=1}^Nc_j\delta_{p_j}$ be a finite signed measure on $\mathbb{R}$,  and consider $U\colon \mathbb{R}^{d}\times \mathbb{R}^d \to \mathbb{R}$ defined by\begin{align*}
    U(x,y)=\int_{\mathbb{R}}\k_d(e^{\frac{s}{2}}x, e^{\frac{s}{2}}y)\mathclose{}\gamma(\textup{d}s)=\int_{\mathbb{R}}\kappa_d(e^s|x||y|)\gamma(\textup{d}s).
\end{align*} Then distributionally, \begin{align*}
    \hat{U}(x,y)=\int_{\mathbb{R}}e^{-sd}\k_d\mathopen{}\left(e^{-\frac{s}{2}}x,e^{-\frac{s}{2}}y\right)\mathclose{}\gamma(\textup{d}s)=\int_{\mathbb{R}}e^{sd}\kappa_d(e^s|x||y|)\gamma^-(\textup{d}s),
\end{align*} where $\gamma^-(A)\coloneqq \gamma(-A)$ is the reflection of $\gamma$ across the origin. So $u\colon \mathbb{R}^d\times \mathbb{R}^d \to \mathbb{R}$ defined by \begin{align*}
    u(x,y)&=\bigintsss_{\mathbb{R}}\kappa_d(e^s|x||y|)\mathopen{}\left(\gamma(\textup{d}s)-e^{sd}\gamma^-(\textup{d}s)\right)\mathclose{} .
\end{align*} is a $-1$-eigendistribution. Let \begin{align*}
    T_d\gamma(\textup{d}s)\coloneqq \gamma(\textup{d}s)-e^{sd}\gamma^-(\textup{d}s).
\end{align*} We then have with this parametrisation, \begin{align*}
    \left\lVert k_d-u\right\rVert_{L^\infty}=\sup_{r>0}\left|\kappa_d(r)-\int_{\mathbb{R}}\kappa_d(re^s)T_d\gamma(\textup{d}s)\right|.
\end{align*} 
Note \begin{align*}
    T_d\gamma(s)&=\sum_{j}c_j\mathopen{}\left(\delta_{p_j}(s)-e^{sd}\delta_{p_j}(-s)\right)\mathclose{}=\sum_{j}c_j\mathopen{}\left(\delta_{p_j}(s)-e^{sd}\delta_{0}(-s-p_j)\right)\mathclose{} \\
    &=\sum_{j}c_j\mathopen{}\left(\delta_{p_j}(s)-e^{sd}\delta_{0}(s+p_j)\right)\mathclose{},
\end{align*} so $ T_d\gamma=\sum_{j}c_j(\delta_{p_j}-e^{-dp_j}\delta_{-p_j})$. This gives\begin{align*}
    \int_{\mathbb{R}}\kappa_d(re^s)T_d\gamma(\textup{d}s)=\sum_{j=1}^Nc_j\mathopen{}\left(\kappa_d(re^{p_j})-e^{-dp_j}\kappa_d(re^{-p_j})\right)\mathclose{},
\end{align*} so we have, for any choice of  $c_j,p_j \in \mathbb{R}$,  \begin{align*}
    w_d^{\rad}\le \sup_{r>0}\left|\sum_{j=1}^Nc_j\mathopen{}\left(\kappa_d(re^{p_j})-e^{-dp_j}\kappa_d(re^{-p_j})\right)\mathclose{}-\kappa_d(r)\right|.
\end{align*} 
For each dimension $d$, for $c,p \in \mathbb{R}^N$, we let \begin{align*}
    \mathcal{I}_d(c,p)\coloneqq \sup_{r>0}\left|\sum_{j=1}^Nc_j\mathopen{}\left(\kappa_d(re^{p_j})-e^{-dp_j}\kappa_d(re^{-p_j})\right)\mathclose{}-\kappa_d(r)\right|,
\end{align*} so we have $w_d^{\rad}\le \mathcal{I}_d(c,p)$. With this formulation, optimising $\mathcal{I}$, Table~\ref{tab:numerics} records (certified) upper bounds for $(w_d^{\rad})^{\frac1{d}}$, with the corresponding lower bounds for $\A_{\sigma}(d)$. The entries $\A_{\sigma}(d)$, $\A_{\sigma}(d)/\sqrt{d}$ were computed from certified upper bounds on $w_d^{\rad}$, and rounded down to $4$ decimal places.  The code is publicly available at this \href{https://github.com/Roni-Edwin/Fourier-numerical-computations/releases/tag/v1.1.1}{GitHub release accompanying this paper}.

 { \small 
	\bibliographystyle{plain}
	\bibliography{bib.bib} }

\clearpage
\begin{table}[p]
\centering
\small
\renewcommand{\arraystretch}{0.95}
\begin{tabular}{|r|r|r|r|}
\hline
$d$ & $(w_d^{\mathrm{rad}})^{1/d}$ & $\A_\sigma(d)$
    & $\A_\sigma(d)/\sqrt d$ \\
    & (upper bound) & (lower bound) & (lower bound) \\
\hline

$4$  & $0.7681$  & $0.6176$ & $0.3088$ \\
           \hline
           $5$ & $0.7674$ & $0.7085$ & $0.3168$ \\
           \hline
           $6$ & $0.7680$ & $0.7860$ & $0.3208$ \\
           \hline
           $7$ & $0.7688$ & $0.8547$ & $0.3230$ \\
           \hline
            $8$  &  $0.7701$ & $0.9165$& $0.3240$ \\
            \hline
             $9$  &  $0.7711$ & $0.9736$& $0.3245$ \\
            \hline
            $10$  &  $0.7713$ & $1.0278$& $0.3250$ \\
            \hline
            $11$ & $0.7721$ & $1.0779$ & $0.3250$ \\
            \hline
              $12$  & $0.7720$ & $1.1266$ & $0.3252$ \\
            \hline
              $13$  & $0.7725$ & $1.1720$ & $0.3250$ \\
            \hline
              $14$  & $0.7728$ & $1.2157$ & $0.3249$ \\
            \hline
              $15$  & $0.7733$ & $1.2572$ & $0.3246$ \\
            \hline
              $16$  & $0.7738$ & $1.2973$ & $0.3243$ \\
            \hline
              $17$  & $0.7745$ & $1.3355$ & $0.3239$ \\
            \hline
              $18$  & $0.7747$ & $1.3732$ & $0.3236$ \\
            \hline
              $19$  & $0.7748$ & $1.4099$ & $0.3234$ \\
            \hline
              $20$  & $0.7752 $  & $1.4452$ & $0.3231$ \\
            \hline
             $22$ & $0.7755$ & $1.5135$ & $0.3226$  \\
\hline
$24$ & $0.7760$ & $1.5783$ & $0.3221$ \\
\hline
$26$ & $0.7764$ & $1.6403$ & $0.3217$ \\
\hline
$28$ & $0.7767$ & $1.6998$ & $0.3212$ \\
\hline
$30$ & $0.7771$ & $1.7571$ & $0.3208$ \\
\hline
$34$ & $0.7777$ & $1.8659$ & $0.3200$ \\
\hline
$38$ & $0.7783$ & $1.9681$ & $0.3192$ \\
\hline
             $42$ & $0.7789$ & $2.0648$ & $0.3186$ \\
\hline
$46$ & $0.7795$ & $2.1565$ & $0.3179$ \\
\hline
$50$ & $0.7801$ & $2.2443$ & $0.3174$ \\
\hline
$54$ & $0.7809$ & $2.3277$ & $0.3167$ \\
\hline
$58$ & $0.7815$ & $2.4082$ & $0.3162$ \\
\hline
$62$ & $0.7826$ & $2.4843$ & $0.3155$ \\
\hline
$66$ & $0.7837$ & $2.5578$ & $0.3148$ \\
\hline
$70$ & $0.7849$ & $2.6282$ & $0.3141$ \\
\hline
$74$ & $0.7867$ & $2.6945$ & $0.3132$ \\
\hline
$78$ & $0.7884$ & $2.7586$ & $0.3123$ \\
\hline
$82$ & $0.7866$ & $2.8335$ & $0.3129$ \\
\hline
$86$ & $0.7864$ & $2.9008$ & $0.3128$ \\
\hline
$90$ & $0.7881$ & $2.9598$ & $0.3119$ \\
\hline
$94$ & $0.7885$ & $3.0219$ & $0.3116$ \\
\hline
$98$ & $0.7884$ & $3.0845$ & $0.3115$ \\
\hline
$100$ & $0.7878$ & $3.1176$ & $0.3117$ \\
\hline

\end{tabular}
\caption{Certified bounds for radial Fourier inequalities and sign uncertainty.}
\label{tab:numerics}
\end{table}
\clearpage

\end{document}